\documentclass{article}

\usepackage{microtype}
\usepackage{graphicx}
\usepackage{subfigure}
\usepackage{booktabs} 

\usepackage{hyperref}

\usepackage[accepted]{icml2021}

\usepackage{hyperref}
\usepackage{url}
\usepackage{enumerate}

\usepackage{amsmath}
\usepackage{amssymb}
\usepackage{amsfonts}
\usepackage{color}
\usepackage{cleveref}

\usepackage{thmtools}
\usepackage{thm-restate}

\usepackage{algorithm}
\usepackage{algorithmic}

\usepackage[normalem]{ulem}

\newcommand{\norm}[1]{\vert \vert {#1} \vert \vert}

\newcommand{\qed}{\(\hfill \square\)}

\newcommand{\p}{\partial}
\newcommand{\inv}[1]{{#1}^{-1}}

\DeclareMathAlphabet{\mathpzc}{OT1}{pzc}{m}{it}
\newcommand{\obj}[1]{{\mathpzc{#1}}}
\newcommand{\pr}[1]{\mathtt{X}_{#1}}

\DeclareMathOperator{\id}{\mathtt{I}}
\DeclareMathOperator{\tw}{tw}
\DeclareMathOperator{\pa}{\operatorname{pa}} 
\DeclareMathOperator{\ch}{\operatorname{ch}}  
\DeclareMathOperator{\nh}{\delta}  

\DeclareMathOperator{\T}{\mathtt{T}}  
\DeclareMathOperator{\iput}{\mathtt{Input}}

\newtheorem{lem}{Lemma}
\newtheorem{theo}{Theorem}
\newtheorem{defn}{Definition}

\usepackage[latin1]{inputenc}
\usepackage{tikz}
\usetikzlibrary{trees, arrows, matrix, positioning}

\icmltitlerunning{Computing the Newton-step faster than Hessian-accumulation}

\begin{document}

\twocolumn[
\icmltitle{Computing the Newton-step faster than Hessian accumulation}



\icmlsetsymbol{equal}{*}

\begin{icmlauthorlist}
  \icmlauthor{Akshay Srinivasan}{sony}
  \icmlauthor{Emanuel Todorov}{uw}
\end{icmlauthorlist}

\icmlaffiliation{sony}{SonyAI, Tokyo, Japan.}
\icmlaffiliation{uw}{University of Washington, Seattle, WA, USA}

\icmlcorrespondingauthor{Akshay Srinivasan}{akssri@vakra.xyz}

\icmlkeywords{Machine Learning, ICML}

\vskip 0.3in
]



\printAffiliationsAndNotice{}  

\begin{abstract}
  Computing the Newton-step of a generic function with $N$ decision variables takes $O(N^{3})$ flops. In this paper, we show that given the computational graph of the function, this bound can be reduced to $O({m \tau}^{3})$, where $\tau, m$ are the width \& size of a tree-decomposition of the graph. The proposed algorithm generalizes nonlinear optimal-control methods based on LQR to general optimization problems and provides non-trivial gains in iteration-complexity even in cases where the Hessian is dense.
\end{abstract}

\section{Introduction}
Newton's method forms the basis for second-order methods in optimization. Computing the Hessian of a generic function \(\obj{f} : \mathbb{R}^N \rightarrow \mathbb{R}\), requires \(O(N^2)\) flops; inverting this matrix requires a further \(O(N^{3})\) flops. This super-linear scaling in compute/memory requirements is a major obstacle in the application of such methods despite their quadratic convergence.

The iteration-complexity of second-order methods has necessitated the development of specialized algorithms for restricted classes of problems. In particular, \emph{Differential Dynamic Programming} (DDP) methods have been proposed for solving nonlinear optimal-control problems, and these methods, quite miraculously, achieve quadratic convergence despite their linear iteration-complexities \cite{jacobson1970differential} \cite{de1988differential} \cite{wright1990solution}.

DDP-methods were historically derived using \emph{Dynamic Programming} (DP) and resemble LQR control-design. They were later, quite surprisingly, also found to be related to both, the Newton-iteration on the unconstrained problem \cite{de1988differential}, and the SQP-iteration on the constrained problem \cite{wright1990solution}. Sadly, the method is not trivially generalized to other objective functions.

Applications of DP to other domains, esp. the solution of least-squares problems and inference in graphical models, have seen great success, but such techniques are not applicable to problems in optimal-control or neural-networks, due to their inability to handle function compositions.

In this paper, we develop a method that exploits the compositional structure in a given objective function in order to automatically derive a \emph{fast} Newton update. This is done by extending the connections of DDP to constrained \& unconstrained optimization using tools from \emph{Automatic Differentiation} (AD), and by using techinques from \emph{Sparse Linear Algebra} (SLA) to bound iteration-complexity.

\section{Problem setup}
\subsection{Computational graph}
Let \(\mathcal{G}\) be a Directed Acyclic Graph (DAG). Define,
\begin{equation}
  \begin{aligned}
    \pa(u) &\triangleq \{v | (v, u) \in E[\mathcal{G}]\}, \quad \mbox{(parents)}\\
    \ch(u) &\triangleq \{v | (u, v) \in E[\mathcal{G}]\}. \quad \mbox{(children)}
  \end{aligned}
\end{equation}

Let every vertex \(v \in V[\mathcal{G}]\) be associated with a \emph{state} \(\pr{v} \in U_v \subset \mathbb{R}^{n_v}\) for some open set $U_v$. Let \(\pr{A}\) be the (labelled) concatenation of \emph{state}s associated with vertices in set \(A \subset V[\mathcal{G}]\).

Let the \emph{input} nodes \(\iput \triangleq \{u_1, u_2,\dots, u_n\} \subset V[\mathcal{G}]\) be the parentless vertices in \(\mathcal{G}\). Let the \emph{states} of the non-input nodes be defined recursively by \(\pr{v} \triangleq \Phi_v(\pr{\pa(v)})\) for some given function \(\Phi_v : \prod_{z \in \pa(v)} U_z \rightarrow U_v\). Since \(\mathcal{G}\) is a DAG, it follows that \(\pr{V[\mathcal{G}]}\) is uniquely determined from the input state \(\pr{\iput}\) and functions \(\{\Phi_v\}_{ v \in V[\mathcal{G}] \backslash X }\).

An objective function \(\obj{f} : U_{x_1} \times \dots U_{x_n} \rightarrow \mathbb{R}\) has the \emph{computational structure} given by the tuple \((\mathcal{G}, \{\Phi_v\}_{ v \in V[\mathcal{G}] \backslash \iput}, \{\obj{l}_v\}_{ v \in V[\mathcal{G}]})\) if it can be written as the sum of local objectives \(\obj{l}_v : \prod_{z \in \{v\} \cup \pa(v)} U_z \rightarrow \mathbb{R}\) on the graph \(\mathcal{G}\) in the following form \eqref{eqn:struct},
\begin{equation}
\label{eqn:struct}
\begin{aligned}
\obj{f}: (\pr{x_1}, \dots, \pr{x_n})& \mapsto \sum_{v \in V[\mathcal{G}]} \obj{l}_v(\pr{v \cup \pa(v)}),\\
\pr{v} \leftarrow \Phi_v(\pr{\pa(v)})&, \quad \forall v \in V[\mathcal{G}], \pa(v) \neq \emptyset.
\end{aligned}
\end{equation}
The \emph{assignment} operator, '\(\leftarrow\)', is explicitly distinguished from the equality operator, '\(=\)', which is taken to represent a constraint in the program. We refer to the DAG \(\mathcal{G}\) as the \emph{computational graph} of \(f(\cdot)\). \\

\textsc{Notation.} The symbolism \(\p_u v \triangleq {\p \pr{v} \over \p \pr{u}} \big|_{\pr{\iput}}\), will be employed for succintly denoting partial derivatives. The derivative operator \emph{w.r.t} the (labelled) set \(A = \{v_1, v_2, \dots\} \subset V[\mathcal{G}]\), will similarly be denoted by \(\p_{A} \triangleq [\p_{a_1}, \p_{a_2}, \dots]\). \\

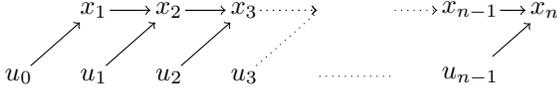
\begin{figure}
  \centering
  \begin{tikzpicture}[shorten >=1pt,->]
    \tikzstyle{vertex}=[circle, minimum size=0pt,inner sep=0pt]

    \foreach \name/\x in {0/0, 1/1, 2/2, 3/3, n-1/6}
    \node[vertex] (U-\name) at (\x,0) {$u_{\name}$};

    \foreach \name/\x in {1/1, 2/2, 3/3, n-1/6, n/7}
    \node[vertex] (X-\name) at (\x, 25pt) {$x_{\name}$};


    \foreach \from/\to in {1/2, 2/3, n-1/n}
    { \draw (X-\from) edge (X-\to); }
    \draw[->, dotted] (X-3) -- (4, 25pt);
    \draw[->, dotted] (X-3) -- (4, 25pt);
    \draw[->, dotted] (5, 25pt) -- (X-n-1);

    \foreach \from/\to in {0/1, 1/2, 2/3, n-1/n}
    { \draw[->] (U-\from) -- (X-\to); }
    \draw[-, dotted] (U-3) -- (4, 25pt);
    \draw[-, dotted] (4, 0) -- (5, 0);
  \end{tikzpicture}
  \caption{(Optimal control) The nodes \(\{x_i\}, \{u_i\}\) represent the states and control inputs of the dynamical system.}
  \label{fig:optcon-directed}
\end{figure}

\subsection{Motivation.} Consider the canonical optimal-control problem,
\begin{equation}
\label{eqn:optcon}
\begin{aligned}
&\min_{u_0, u_1, \dots, u_{n - 1}} \left[\obj{J}(u_0, \dots, u_n) \triangleq \sum_{i = 0}^{n - 1} \obj{l}_i(x_i, u_i) + \obj{l}_n(x_n)\right],\\
&\quad \forall i, x_{i + 1} \leftarrow \mathtt{f}(x_i, u_i),
\end{aligned}
\end{equation}
wherein, \(\mathtt{f}(\cdot, \cdot)\) denotes the system dynamics, and \(\obj{l}_i(\cdot, \cdot)\)'s denote the losses incurred. The order of computation for the objective \eqref{eqn:optcon} can be represented by a linear-chain as shown in \Cref{fig:optcon-directed}.

The computational graph \Cref{fig:optcon-directed}, while sparse, does not imply sparsity in the objective. Substitution of the symbol assignments in the unconstrained objective, has a cascading effect whereby variables get \emph{transport}ed down \(\mathcal{G}\),
\begin{flalign*}
\obj{J}(u_0, \dots, u_n) = &\obj{l}_0(x_0, u_0) +\\
&\obj{l}_1(\mathtt{f}(x_0, u_0), u_1) +\\
&\obj{l}_2(\mathtt{f}(\mathtt{f}(x_0, u_0), u_1), u_2) + \dots.
\end{flalign*}
This results in a Hessian that is both fully dense and (very) badly conditioned. These two properties render both iterative and direct-factorization methods unsuitable for solving this problem.

The computational graph does, however, fully encode the sparsity of the Lagrangian of the constrained program; the constrained form being obtained by replacing '\(\leftarrow\)' with '\(=\)', and including \(\{x_i\}'s\) in the domain of optimization. This results in a \emph{fast} linear-time SQP (Sequential Quadratic Programming) iteration for the optimal-control problem but increases implementation complexity by requiring \emph{dual} and \emph{non-input} updates.

This presents one with a strange set of choices: slow-simple unconstrained optimization or fast-complex constrained optimization. DDP-like methods resolve this dichotomy for optimal-control problems \cite{jacobson1970differential} \cite{de1988differential} \cite{wright1990solution}. Our goal is to do the same for general problems.

\section{Graphical Newton}
\label{sec:newt}
Consider the objective function in \eqref{eqn:struct}, defined by the tuple \((\mathcal{G}, \{\Phi_v\}, \{\obj{l}_v\})\). The optimization problem of interest is the following,
\begin{equation}
  \label{eqn:obj}
\begin{aligned}
\min_{\{\pr{v} | \forall v \in \iput\}} &\left( \obj{f} \triangleq \sum_{v \in V[\mathcal{G}]} \obj{l}_v(\pr{v \cup \pa(v)}) \right),\\
\pr{v} \leftarrow \Phi_v(\pr{\pa(v)}),& \quad \forall v \in V[\mathcal{G}], \pa(v) \neq \emptyset.
\end{aligned}
\end{equation}
The corresponding constrained problem is obtained by simply replacing the operator '\(\leftarrow\)' with '\(=\)' in \eqref{eqn:obj} and enlarging the domain of optimization.
\begin{equation}
  \label{eqn:coobj}
\begin{aligned}
\min_{\{\pr{v} | \forall v \in V[\mathcal{G}]\}} &\left( \obj{f} \triangleq \sum_{v \in V[\mathcal{G}]} \obj{l}_v(\pr{v \cup \pa(v)}) \right),\\
\pr{v} = \Phi_v(\pr{\pa(v)}),& \quad \forall v \in V[\mathcal{G}], \pa(v) \neq \emptyset.
\end{aligned}
\end{equation}

The constrained program (\ref{eqn:coobj}) induces the following Lagrangian function,
\begin{equation}
  \label{eqn:cdef}
  \begin{aligned}
    &\mathcal{L}(\pr{V[\mathcal{G}]}, \lambda_{V[\mathcal{G}] \backslash \iput}) \triangleq
    \sum_{v \in V[\mathcal{G}]} \obj{l}_v(\pr{v \cup \pa(v)}) + \sum_{\substack{v \in V[\mathcal{G}],\\ \pa(v) \neq \emptyset}} \lambda_v^{\mathtt{T}} \mathbf{h}_v(\pr{v \cup \pa(v)}),\\
    &\mbox{where,} \; \mathbf{h}_v(\pr{v \cup \pa(v)}) \triangleq \Phi_v(\pr{\pa(v)}) - \pr{v}, \; \forall v \in V[\mathcal{G}], \pa(v) \neq \emptyset.
  \end{aligned}
\end{equation}
The necessary first-order conditions for optimality are given by,
\begin{equation}
  \label{eqn:lagrangefirst}
  \p_{V} \mathcal{L}(\pr{V}^*, \lambda_V^*) = 0,\quad h(\pr{V}^*) = 0.
\end{equation}

Linearization of the first-order conditions for this constrained problem (around a nominal \((\tilde{\pr{V}}, \tilde{\lambda})\)) yields a system of KKT equations, whose solution yields the SQP search direction.
\begin{equation}
  \label{eqn:csqpstep}
  \left[\begin{array}{c c}
          \p^2_{V} \mathcal{L} & \p_{V} \mathbf{h}^{\mathtt{T}} \\
          \p_{V} \mathbf{h} & 0
    \end{array}\right] \left[\begin{array}{c} \delta \pr{V} \\ \lambda^{+}
    \end{array}\right] = \left[\begin{array}{c} - \p_{V} \obj{f} \\ -\mathbf{h}
                               \end{array} \right],
\end{equation}
where \(\lambda^+ \triangleq \tilde{\lambda} + \delta \lambda\). The sequence of iterates obtained by taking appropriate steps along \((\delta \pr{V}, \delta \lambda)\), converges quadratically near a strongly-convex local minimum.

The principal result of this paper is the connection between the unconstrained and constrained formulations,

\begin{restatable}[Newton direction]{theo}{gnewt}
  \label{thm:gnewt} The Newton direction for unconstrained problem \eqref{eqn:obj} is given by an iteration of SQP for the constrained problem \eqref{eqn:coobj}, when \(\pr{V[\mathcal{G}]}\) is feasible and when \(\forall v, \lambda_v = \p_v \obj{f}\) \eqref{eqn:a_rmad}.
\end{restatable}
\emph{Proof.} See \Cref{sec:gnewtproof}.

\Cref{thm:gnewt} can trivially be extended to arbitrary objective functions on the graph $\mathcal{G}$ and is not restricted to form \eqref{eqn:coobj} where objective functions have the same sparsity as the underlying computation. The dual values $\lambda_v = \p_v \obj{f}$ can be computed in linear-time with reverse-mode AD.

\Cref{thm:gnewt} immediately leads to the sparse Newton-iteration given in \Cref{alg:gnewt}.

\begin{algorithm}
  \caption{Graphical Newton}
  \label{alg:gnewt}
\begin{algorithmic}[1]
  \STATE {\bfseries Input:} initial $\pr{\iput}^0$, tuple $(\mathcal{G}, \{\Phi_v\}, \{\obj{l}_v\})$
  \REPEAT
  \STATE Compute non-inputs, local objectives, and their derivatives (\cref{eqn:obj}).
  \STATE Set dual-values to \(\lambda_v = \p_v \obj{f}, \forall v\) by reverse-mode AD (\cref{eqn:a_rmad}).
  \STATE Solve the KKT-system (\cref{eqn:csqpstep}).
  \STATE Compute step-length \(\eta\) \emph{via} linesearch on inputs \(\pr{\iput}\) (only).
  \STATE Update inputs (only): $\pr{\iput} \leftarrow \pr{\iput} + \eta \delta \pr{\iput}$
  \UNTIL{$\norm{\p_{\iput} \obj{f}} \leq \epsilon$ }
\end{algorithmic}
\end{algorithm}

\Cref{alg:gnewt} differs markedly from both SQP and current techniques in AD. The non-input primals are computed by running a forward pass on the computation graph, while the dual-values are set to fixed values by running reverse-mode AD. Once the KKT system is solved, the algorithm proceeds to update the primal inputs only, without requiring non-input updates, dual updates, contraint penalties, or any of the other machinery from constrained optimization.

The efficiency of \Cref{alg:gnewt} stems from the fact that the KKT matrix is typically much sparser than the Hessian. It can be further shown that the common Hessian-accumulation/inversion method is equivalent to a particular, generally a (very) suboptimal, pivot-ordering for solving the sparse KKT system.\\

\subsection{KKT complexity}
The run-time of the \Cref{alg:gnewt} depends crucially on the time taken to solve the sparse KKT system \eqref{eqn:csqpstep} at every iteration. The complexity of solving such sparse systems depends inturn on the support-graph of the underlying matrix. 

It has been observed that many problems defined locally on graphs, can be solved in linear-time on trees using dynamic programming. These techniques can be extended to general graphs by grouping vertices in such a way as to \emph{mimic} a tree. The size of the largest cluster in this \emph{tree-decomposition} - termed \emph{tree-width} - governs the dominant term in the time-complexity of the resulting overall algorithm. Computing the tree-decomposition with minimal tree-width is NP-hard, but heuristics for finding elimination orderings are often known to do well in practice. 

Run-time complexities in terms of the tree-width are well-established for closely related problems such as Cholesky decomposition, but they appear to be unknown for structured KKT systems such as \eqref{eqn:csqpstep}. The following theorem establishes the required bound for structured-KKT systems arising in \Cref{alg:gnewt}.
\begin{restatable}[KKT tree-width]{theo}{twlem}
\label{thm:twlem} The KKT system \eqref{eqn:csqpstep} associated with the constrained problem \eqref{eqn:coobj} can be solved in time \(O(m \tw(\mathcal{G})^3)\), given the tree-decomposition.
\end{restatable}
\emph{Proof.} See \Cref{sec:twlemproof}.

\section{Special cases}
\subsection{Optimal control}
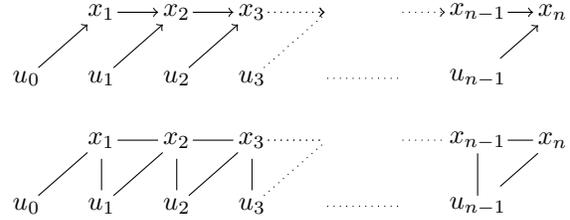
\begin{figure}
  \centering
  \begin{tikzpicture}[shorten >=1pt,->]
    \tikzstyle{vertex}=[circle, minimum size=0pt,inner sep=0pt]

    \foreach \name/\x in {0/0, 1/1, 2/2, 3/3, n-1/6}
    \node[vertex] (U-\name) at (\x,0) {$u_{\name}$};

    \foreach \name/\x in {1/1, 2/2, 3/3, n-1/6, n/7}
    \node[vertex] (X-\name) at (\x, 25pt) {$x_{\name}$};


    \foreach \from/\to in {1/2, 2/3, n-1/n}
    { \draw (X-\from) edge (X-\to); }
    \draw[->, dotted] (X-3) -- (4, 25pt);
    \draw[->, dotted] (X-3) -- (4, 25pt);
    \draw[->, dotted] (5, 25pt) -- (X-n-1);

    \foreach \from/\to in {0/1, 1/2, 2/3, n-1/n}
    { \draw[->] (U-\from) -- (X-\to); }
    \draw[-, dotted] (U-3) -- (4, 25pt);
    \draw[-, dotted] (4, 0) -- (5, 0);
  \end{tikzpicture}

    \begin{tikzpicture}[shorten >=1pt,-]
    \tikzstyle{vertex}=[circle, minimum size=0pt,inner sep=0pt]

    \foreach \name/\x in {0/0, 1/1, 2/2, 3/3, n-1/6}
    \node[vertex] (U-\name) at (\x,0) {$u_{\name}$};

    \foreach \name/\x in {1/1, 2/2, 3/3, n-1/6, n/7}
    \node[vertex] (X-\name) at (\x, 25pt) {$x_{\name}$};


    \foreach \from/\to in {1/2, 2/3, n-1/n}
    { \draw (X-\from) edge (X-\to); }
    \draw[-, dotted] (X-3) -- (4, 25pt);
    \draw[-, dotted] (X-3) -- (4, 25pt);
    \draw[-, dotted] (5, 25pt) -- (X-n-1);

    \foreach \from/\to in {0/1, 1/2, 2/3, n-1/n}
    { \draw[-] (U-\from) -- (X-\to); }
    \foreach \from/\to in {1/1, 2/2, 3/3}
    { \draw[-] (U-\from) -- (X-\to); }
    \draw[-] (6, 4pt) -- (6, 21pt);

    \draw[-, dotted] (U-3) -- (4, 25pt);
    \draw[-, dotted] (4, 0) -- (5, 0);
  \end{tikzpicture}

  \caption{Structure of the optimal control problem as defined in \eqref{eqn:optcon_c}. Top: Computational graph. Bottom: Constraint graph.}
  \label{fig:optcon}
\end{figure}

Consider, again, the canonical optimal control problem from \eqref{eqn:optcon},
\begin{equation}
\label{eqn:optcon_c}
\begin{aligned}
&\min_{u_0, u_1, \dots, u_{n - 1}} \left[\obj{J}(u_0, \dots, u_{n-1}) \triangleq \sum_{i = 0}^{n - 1} \obj{l}_i(x_i, u_i) + \obj{l}_n(x_n)\right],\\
&\quad \forall i, x_{i + 1} \leftarrow \mathtt{f}(x_i, u_i).
\end{aligned}
\end{equation}
The computational graph and its moralized relative for this problem are shown in \Cref{fig:optcon}. The constraint graph is chordal, and permits multiple optimal elimination orderings.

\textsc{Dynamic Programming.}
\begin{table*}[ht]
  \caption{DDP methods and their key differences.}
  \center
  \begin{tabular}{| c | c | c |}
    \hline
    Method & \(\lambda_i\) & \emph{back-substitution} \\
    \hline
    DDP \cite{jacobson1970differential} & \(\p_{x_{i + 1}} V_{i + 1} \cdot \p_{x_i}{\mathtt{f}} \) & Non-Linear \eqref{eqn:nonlinbs}\\
    Stagewise-Newton \cite{de1988differential} & \(\p_{x_i} \obj{J}\) & Linear\\
    Nonlinear Stagewise-Newton \cite{liao1992advantages} & \(\p_{x_i} \obj{J}\) & Non-Linear \eqref{eqn:nonlinbs}\\
    iLQR/iLQG \cite{li2004iterative} & 0 & Non-Linear \eqref{eqn:nonlinbs}\\
    \hline
  \end{tabular}
  \label{tbl:ddpzoo}
\end{table*}

Let's consider the LQR order, \(x_{n}, u_{n-1}, x_{n-1}\dots, u_{0}\). The principal minor of the KKT matrix, corresponding to the clique \(\{x_{n}, x_{n-1}, u_{n-1}, \lambda_{n - 1}^+\}\) is given by,
\begin{equation}
  \label{eqn:minor}
  \begin{aligned}
    &\begin{array}{l c c c | c}
      \p^2_x V_n &  &  & - \id^{\T} & - \p_x V_n \\
     & \p_x^2 \mathcal{H}_{n -1}  & \p_{xu}^2 \mathcal{H}_{n -1} & \p_x \mathtt{f} & - \p_x \obj{l}_{n - 1}\\
     & \p_{ux}^2 \mathcal{H}_{n -1} & \p_{uu}^2 \mathcal{H}_{n -1} & \p_u \mathtt{f} & - \p_u \obj{l}_{n - 1}\\
    -\id^{\T} & \p_x \mathtt{f} & \p_u \mathtt{f} & 0 & 0
    \end{array},\\
    &\mbox{where,} V_n = l_n, \mathcal{H}_{n-1}(x, u) = \obj{l}_{n-1}(x, u) + \lambda_{n-1} \cdot \mathtt{f}(x, u).
  \end{aligned}
\end{equation}
Note that we can only eliminate variables \(\{x_{n}, u_{n-1}\}\) in the above, since \(x_{n-1}\) is part of the separator \emph{i.e} it is connected to nodes outside the clique.\\
Eliminating \(\delta x_{n} =  \inv{\p^2 V_{n}} ( -\p_x V_{n} + \lambda_{n -1 }^+)\),
  \[\begin{array}{c c c | c}
      \p_x^2 \mathcal{H}_{n -1}  & \p_{xu}^2 \mathcal{H}_{n -1} & \p_x \mathtt{f}^{\T} & - \p_x \obj{l}_{n - 1}\\
      \p_{ux}^2 \mathcal{H}_{n -1} & \p_{uu}^2 \mathcal{H}_{n -1} & \p_u \mathtt{f}^{\T} & - \p_u \obj{l}_{n - 1}\\
      \p_x \mathtt{f} & \p_u \mathtt{f} & -\inv{\p^2 V_{n}} & - \inv{\p^2 V_{n}} \p V_{n}
    \end{array},\]
  Eliminating \(\lambda^{+}_{n - 1} = - \p^2V_{n} (- \inv{\p^2 V_{n}} \p V_{n} - \p_x \mathtt{f} \delta x_{n - 1} - \p_u \mathtt{f}_u \delta u_{n-1})\), %
  \[\begin{array}{c c | c}
      Q_{xx} & Q_{ux}^{\T} & - q_x\\
      Q_{ux} & Q_{uu} & - q_u
    \end{array},\]
  where,
  \[\begin{array}{l} Q_{ab} = \p_{ab}^2 \mathcal{H}_{n -1}  + \p_{a} \mathtt{f}\; \p^2 V_{n}\; \p_b \mathtt{f}^{\T},\\
      q_{a} = \p_a \obj{l}_{n - 1} + \p_a \mathtt{f}^{\T} \; \p V_{n}.
    \end{array}\]
  
Eliminate \(\delta u_{n-1} = \inv{Q_{uu}} (-Q_u - Q_{u x} \delta x_{n-1}) \),
\[\begin{array}{c | c}
    Q_{xx} + Q_{ux}^{\T} \inv{Q_{uu}} Q_{ux} & - (Q_x + Q_{ux}^{T} \inv{Q_{uu}} Q_u)
  \end{array},\]
Re-write the above as,
\[\begin{array}{c | c}
    \p_x^2 V_{n - 1} & - \p_x^2 V_{n-1}.
  \end{array},\]
The eliminations of the adjoining clique can be carried out in a similar manner.

This procedure is identical to the \emph{backward pass} in DDP-methods.\footnote{These computations can be carried out more efficiently \& stably in the square-root form.} In the backsubstitution phase, also known as \emph{forward pass}, the non-linearities can be used directly without having to resort to the use of linearizations,
\begin{equation}
  \label{eqn:nonlinbs}
  \begin{aligned}
    \delta x_{n} &=  \inv{\p^2 V_{n}} ( -\p_x V_{n} + \lambda_{n -1 }^+)\\
    &= \inv{\p^2 V_{n}} ( -\p_x V_{n} - \p^2V_{n} (- \inv{\p^2 V_{n}} \p V_{n} \\
    & \hspace{6em}- \p_x \mathtt{f} \delta x_{n - 1} - \p_u \mathtt{f}_u \delta u_{n-1}))\\
    &= \p_x \mathtt{f} \delta x_{n - 1} + \p_u \mathtt{f}_u \delta u_{n-1}.\\
    &\approx \mathtt{f}(x_{n-1} + \delta x_{n - 1}, u_{n -1} + \delta u_{n-1}) - \mathtt{f}(x_{n-1}, u_{n -1}).
  \end{aligned}
\end{equation}

The \Cref{tbl:ddpzoo} illustrates how variations in \emph{forward} \& \emph{backward} passes, correspond to known DP algorithms. The theory also makes it trivial to generalize these algorithms to higher-order dynamics and constrained problems \cite{srinivasan2015}.

\section{Discussion}
The method presented in this paper can be used to compute the Newton step in time \(O(m \tw^3)\), where \(\tw, m,\) are the width and size of a tree-decomposition of the computation graph. The method generalizes many specialized DDP algorithms in numerical optimal-control.

However, while the work presented here provides, in a sense, optimal iteration-complexities, many real-world problems in machine-learning also have large tree-widths. The KKT-matrix factorization of such problems is also quite rife with redundant computation, and indicates the necessity for partial symbolic-condensation in sparse LDL (and QP) solvers. These are topics for future work.



\bibliography{references}
\bibliographystyle{icml2021}

\newpage 
\appendix
  \section{Graphical Newton}
  \label{sec:gnewtproof}
  This section presents a proof of \cref{thm:gnewt}. We proceed by first recalling the $1^{st}, 2^{nd}$-order relations from AD for the objective \eqref{eqn:obj}; this is then related to the problem of computing the solutiong to the KKT system in \eqref{eqn:csqpstep} thus completing the proof.

\subsection{Derivative relations on $\mathcal{G}$}

\textsc{Reverse-mode AD.} The first derivative of the objective \(\obj{f}(\cdot)\) \eqref{eqn:obj} can be calculated by applying the chain rule over \(\mathcal{G}\),
\begin{equation}
  \label{eqn:a_rmad}
  \begin{aligned}
    &\forall v,\quad  \p_v \obj{f} = \sum_{s \in v \cup \ch(v)} \p_v \obj{l}_s + \sum_{d \in \ch(v)} \p_d \obj{f}^{\mathtt{T}} \; \p_v d;\\
    &\quad\quad\quad v \in \pa(d) \Rightarrow \p_v d \triangleq {\p \Phi_d(\pr{\pa(d)}) \over \p \pr{v}}.
  \end{aligned}
\end{equation}
Since \(\mathcal{G}\) is a DAG, there exist child-less nodes (\emph{i.e} \(\ch(v) = \emptyset\)) from which the recursion can be initialized. The recursion then proceeds backward on \(\mathcal{G}\) in a breadth-first order. This algorithm is known in literature as reverse-mode AD.

\textsc{Hessian-vector AD}: For a given infinitesimal change \(\delta {\pr{\iput}}\) in the inputs, the first derivatives exhibit a first-order change \(\delta [\p_v \obj{f}] \triangleq \p_{\iput, v}^2 \obj{f} \cdot \delta \pr{\iput}\), given by the Hessian-vector product. Applying chain-rule over \(\mathcal{G}\) again for all terms in \eqref{eqn:a_rmad} we obtain,
\begin{equation}
  \label{eqn:a_armad}
  \begin{aligned}
    & \forall v,\quad  \delta [\p_v \obj{f}] = \sum_{s \in v \cup \ch(v)} \left(\sum_{a \in v \cup \pa(s)} \p_{va}^2 \obj{l}_s \cdot \delta \pr{a} \right) +\\
    & \sum_{d \in \ch(v)} \left(\delta [\p_d \obj{f}]^{\mathtt{T}} \p_v d + \sum_{a \in \pa(d)} (\p_d \obj{f}^{\mathtt{T}} \; \p_{va}^2 d) \cdot \delta \pr{a})\right);\\
    &\mbox{where,} \\
    &\forall a,\quad \delta \pr{a} = \sum_{d \in \pa(a)} \p_{d} a \cdot \delta \pr{{d}}.
  \end{aligned}
\end{equation}
These equations can be solved, for a given \(\delta \pr{\iput}\) by a forward-backward recursion similar to the one used for solving \eqref{eqn:a_rmad}. Computing the Hessian-vector product in this manner takes time \(O(|E[\mathcal{\hat{G}}]|\omega(\mathcal{\hat{G}})^2)\), where \(\omega(\mathcal{\hat{G}})\) is the clique number of the moralization of \(\mathcal{G}\).

\subsection{Newton direction}
Computing the Newton step requires inverting the Hessian-vector AD process: find \(\delta {\pr{\iput}}\) such that, \(\delta [\p_{\iput} \obj{f}] = - \p_{\iput} \obj{f}\). The inversion of these relations is related to the computation of the SQP direction via \Cref{thm:gnewt},

\gnewt*

\emph{Proof.} The second equation in \eqref{eqn:a_armad} is equivalent to \(\p_V \Phi \cdot \delta S_V = -\Phi\), in \eqref{eqn:csqpstep}. Rearranging the first equation from \eqref{eqn:a_armad}, and setting \(\delta[\p_v \obj{f}] = -\p_v \obj{f}\) for all inputs, we obtain \(\forall v\),
\begin{equation}
\begin{aligned}
\label{eqn:a_ppobj}
0 =&\sum_{\substack{s \in v \cup \ch(v),\\ a \in v \cup \pa(s)}} \p_{va}^2 \obj{l}_s \;\delta \pr{a} + \sum_{\substack{d \in \ch(v),\\ a \in \pa(d)}} (\p_d \obj{f}^{\mathtt{T}}\; \p_{va}^2 d) \; \delta \pr{a} +\\
& - \left(\begin{cases} \delta [\p_v \obj{f}] & \pa(v) \neq \emptyset\\ - \p_v \obj{f} & \mbox{otherwise}\end{cases}\right) + \sum_{d \in \ch(v)}  (\p_v d)^{\mathtt{T}} \delta [\p_d \obj{f}].
\end{aligned}
\end{equation}
Similarly, expanding the top block in \eqref{eqn:csqpstep} using the definitions in \eqref{eqn:obj} \& \eqref{eqn:csqpstep}, we obtain \(\forall v\),
\begin{equation}
\begin{aligned}
\label{eqn:a_ppLexp}
-\p_{v} \mathcal{L} =&\sum_{\substack{s \in v \cup \ch(v),\\a \in v \cup \pa(s)}} \p_{va}^2 \obj{l}_s \; \delta \pr{a} + \sum_{\substack{d \in \ch(v),\\a \in \pa(d)}} (\lambda_d^{\mathtt{T}} \; \p_{av}^2 d) \; \delta \pr{a} +\\
& - \left(\begin{cases} \delta \lambda_v & \pa(v) \neq \emptyset \\ 0 & \mbox{otherwise}\end{cases} \right) + \sum_{d \in \ch(v)}  (\p_v d)^{\mathtt{T}} \; \delta \lambda_d ,
\end{aligned}
\end{equation}
where,
\begin{equation}
\begin{aligned}
  \label{eqn:a_pLexp}
  \p_{v} \mathcal{L} = &\begin{cases} \sum_{s \in v \cup \ch(v)}  \p_v \obj{l}_s + \sum_{d \in \ch(v)} \lambda_d \cdot \p_{v} d, & \pa(v) = \emptyset\\ \begin{array}{l} \sum_{s \in v \cup \ch(v)}  \p_v \obj{l}_s + \sum_{d \in \ch(v)} \lambda_d \cdot \p_{v} d \\- \lambda_v \end{array}, & \mbox{otherwise} \\  \end{cases}
\end{aligned}
\end{equation}
The result follows from equations \eqref{eqn:a_rmad}, \eqref{eqn:a_ppobj}, \eqref{eqn:a_ppLexp} \& (\cref{eqn:a_pLexp}).

\qed

\section{KKT complexity}
\label{sec:twlemproof}




In this section, we provide a Message Passing [MP] algorithm for solving KKT systems arising in \Cref{alg:gnewt} and show that it has a theoretical run-time bound of \(O(m \tw^3)\), given a tree-decomposition. %



\subsection{Hypergraph structured QPs}
For a hypergraph \(\mathcal{H}\), denote the adjacency and incidence matrices by \(\mathcal{A}[\mathcal{H}]\) \& \(\mathcal{B}[\mathcal{H}]\) respectively,
\begin{equation}
\begin{aligned}
  &\mathcal{A}[\mathcal{H}] \in \mathbb{R}^{|V[\mathcal{H}]|\times|V[\mathcal{H}]|}, \quad\quad \mathcal{B}[\mathcal{H}] \in \mathbb{R}^{|E[\mathcal{H}]|\times|V[\mathcal{H}]|},\\
  &\mathcal{A}[\mathcal{H}]_{uv} = \begin{cases} 1 & \exists e \in E[\mathcal{H}],\; u, v \in e\\ 0 & \mbox{otherwise}\end{cases}\\
  &\mathcal{B}[\mathcal{H}]_{eu} = \begin{cases} 1 & u \in e\\ 0 & \mbox{otherwise}\end{cases}
\end{aligned}
\end{equation}
Given such a hypergraph \(\mathcal{H}\), the family of QPs we're interested in solving is the following,
\begin{equation}
  \label{eqn:spqpfamily}
  \begin{aligned}
    &\min_{x} \sum_{e \in E[\mathcal{H}]}{1 \over 2} \pr{e}^{\mathtt{T}} Q_e \pr{e} - b_e^{\mathtt{T}} \pr{e}, \\
    &\forall e \in E[\mathcal{H}], \quad G_{e} \pr{e} = h_e.
  \end{aligned}
\end{equation}
Assuming that the QP has a bounded solution and that the constraints are full rank, the minimizer to (\cref{eqn:spqpfamily}) is given by the solution to the following KKT system,
\begin{equation}
  \label{eqn:spfamily}
  \begin{aligned}
    \left[\begin{array}{c c} Q & G^{\mathtt{T}} \\ G & 0\end{array}\right] \left[\begin{array}{c} x \\ \lambda\end{array}\right] &=
    \left[\begin{array}{c} b \\ h \end{array}\right],\\
    x, b \in \mathbb{R}^{|V|}, \quad& \lambda, h \in \mathbb{R}^{M},\\
\end{aligned}
\end{equation}
where \(Q, G, \lambda, x, b\) are concatenation of terms defined in \eqref{eqn:spqpfamily} respectively. The sparsity/support of (\cref{eqn:spfamily}) is closely related to \(\mathcal{H}\) because the quadratic part of the KKT equation has the sparsity of the adjacency matrix, and the row of the constraint \(G_{i, :}\) has the same sparsity as some edge \(e \in E[\mathcal{H}]\).

\textsc{Tree decomposition}:
Extending the notion of Dynamic Programming to non-trees (including Hypergraphs) requires a partitioning of the graph so as to satisfy a \emph{lifted} notion of being a tree. Tree decomposition captures the essence of such graph partitions,

\begin{defn}{(Tree decomposition)}
A tree-decomposition of a hypergraph \(\mathcal{H}\) consists of a tree \(\mathcal{T}\) and a map \(\chi : V[\mathcal{T}] \rightarrow 2^{V[\mathcal{H}]}\), such that,
\begin{enumerate}[i]
\item (\emph{Vertex cover}) \(\cup_{i \in V[\mathcal{T}]} \chi(i) = V[\mathcal{H}].\)
\item (\emph{Edge cover}) \(\forall e \in E[\mathcal{H}],\; \exists i \in V[\mathcal{T}], e \subset \chi(i).\)
\item (\emph{Induced sub-tree}) \(\forall u \in V[\mathcal{H}],\; \mathcal{T}_u \triangleq \mathcal{T}[\{i \in V[\mathcal{T}]| u \in \chi(i)\}]\; \mbox{is a non-empty subtree}\)
\end{enumerate}
The tree-width of a tree-decomposition \(\mathcal{T}\) is defined to be \(\tw(\mathcal{T}) = \max_{v \in V[\mathcal{T}]} |\chi(v)| - 1\). The tree-width of a graph \(\mathcal{H}\) is defined to be the minimal tree-width attained by any tree-decomposition of \(\mathcal{H}\).
\end{defn}

We define the vertex-induced subgraph in what follows to be \(\mathcal{H}[S] \triangleq (V[\mathcal{H}], \{e \cap S, e\in E[\mathcal{H}]\})\). The following lemma ensures that such a decomposition ensures \emph{local dependence}.
\begin{lem}
\label{lem:esep}{(Edge separation)} Deleting the edge \(xy \in E[\mathcal{T}]\), renders \(\mathcal{H}[V \backslash (\chi(x) \cap \chi(y))]\) disconnected.
\end{lem}

\textsc{Hypertree structured QP}: The tree-decomposition itself can be considered a Hypergraph, \((V[\mathcal{H}], \{\chi(u), \forall u \in V[\mathcal{T}]\})\). Such a \emph{Hypertree}\footnote{There are multiple definitions of a \emph{Hypertree}; we use the term to mean a maximal Hypergraph, whose tree-decomposition can be expressed in terms of its edges.} can also be thought of as a chordal graph. We assume henceforth that the given graph \(\mathcal{H}\) is a hypertree, and that \(\mathcal{T}\) is its tree-decomposition.

A message passing [MP] algorithm for solving \eqref{eqn:spfamily} on such Hypertrees is given below. The gather stage of the message passing algorithm is illustrated in \cref{alg:gqp} \footnote{Note that the addition is performed vertex label-wise in Line~6 of Algorithm~\cref{alg:gqp}.}. The function, Factorize, computes the partial LU decomposition of its arguments; we describe below, its operation.

Denote the vertices that are interior to \(l\) by \(\iota = \chi(l) \cap \chi(p)\), and those on the boundary (\emph{i.e} common to \(p, l\)) by \(\xi = \chi(l) \backslash \chi(p)\), and let \(r = \operatorname{rank}(\tilde{Q}_{\iota, \iota})\).  The function computes Gaussian-BP messages from block pivots \(\{2, 3\}\) to \(\{1, 4\}\) in \eqref{eqn:mpcore}. Note that, unlike Gaussian-BP, the matrices in \eqref{eqn:mpcore} are not necessarily positive definite, but are however invertible.

\begin{algorithm}[h]
  \caption{Graphical QP}
  \label{alg:gqp}
\begin{algorithmic}[1]
  \STATE {\bfseries Given}: $\mathcal{T}, \mathcal{H}, \{Q_e\}, \{b_e\}, \{G_e\}, \{h_e\}.$
  \STATE
  \STATE {\bfseries function} GatherMessage($l, p, \mathcal{T}$)
  \STATE $(\tilde{Q}_l, \tilde{b}_l, \tilde{G}_l, \tilde{h}_l) \leftarrow (Q_l, b_l, G_l, h_l)$
  \FOR {$c \in \nh_{\mathcal{T}}(l)\backslash p$}
  \STATE $(Q_{c\rightarrow l}, G_{c\rightarrow l}, b_{c\rightarrow l}, h_{c\rightarrow l}) \leftarrow \mbox{GatherMessage}(c, p, \mathcal{T})$
  \STATE $(\tilde{Q_l}, \tilde{b_l}) \leftarrow (\tilde{Q_l}, \tilde{b_l}) + (Q_{c\rightarrow l}, b_{c\rightarrow l})$
  \STATE $\tilde{G}_l \leftarrow [\tilde{G}_l; G_{c\rightarrow l}], \tilde{h}_l \leftarrow [\tilde{h}_l; h_{c\rightarrow l}]$
  \ENDFOR
  \STATE {\bfseries return} Factorize($\chi(l), \chi(p), \tilde{Q}_l, \tilde{b}_l, \tilde{G}_l, \tilde{h}_l$)
  \STATE
  \STATE {\bfseries function} Factorize($\chi(l), \chi(p), \tilde{Q}, \tilde{b}, \tilde{G}, \tilde{h}$)
  \STATE $(\xi, \iota) \leftarrow (\chi(l) \backslash \chi(p), \chi(l) \cap \chi(p))$
  \STATE $r \leftarrow \operatorname{rank}(\tilde{Q}_{\iota, \iota})$
  \STATE {\bfseries return} Gaussian-BP messages from \eqref{eqn:mpcore}.
  \STATE {\bfseries return}
\end{algorithmic}
\end{algorithm}

\begin{equation}
  \label{eqn:mpcore}
\begin{tikzpicture}
  \matrix [matrix of math nodes,left delimiter={[},right delimiter={]}] (m)
  {
    \tilde{Q}_{\xi \xi} & \tilde{Q}_{\iota \xi}^{\mathtt{T}} & \tilde{G}_{:r, \xi}^{\mathtt{T}} & \tilde{G}_{r:, \iota}^{\mathtt{T}}\\
    \tilde{Q}_{\iota \xi} & \tilde{Q}_{\iota \iota} & \tilde{G}_{:r, \iota}^{\mathtt{T}}& \tilde{G}_{r:, \iota}^{\mathtt{T}}\\
    \tilde{G}_{:r, \xi} & \tilde{G}_{:r, \iota} & 0 & 0\\
    \tilde{G}_{r:, \xi} & \tilde{G}_{r:, \iota} & 0 & 0\\
  };
  \draw[color=red] (m-2-2.north west) -- (m-2-3.north east) -- +(0, -1.3) -- +(-1.75, -1.3) -- (m-2-2.north west);
  \matrix [matrix of math nodes,left delimiter={[},right delimiter={]}] (v) at (3, 0)
  {
    \pr{\xi}\\
    \pr{\iota}\\
    \lambda_{:r}\\
    \lambda_{r:}\\
  };
  \node at (4, 0) {$=$};
  \matrix [matrix of math nodes,left delimiter={[},right delimiter={]}] (v) at (5, 0)
  {
    \tilde{b}_{\xi}\\
    \tilde{b}_{\iota}\\
    \tilde{h}_{:r}\\
    \tilde{h}_{r:}\\
  };
\end{tikzpicture}
\end{equation}



Gaussian Belief-Propagation is essentially a re-statement of LU decomposition, 
and consists of messages of the form, 
\begin{equation}
  \begin{aligned}
    \mu_{i \rightarrow j} &:= [J_{i \rightarrow j}, h_{i \rightarrow j}] = [J_{ii}, h_i] - \sum_{k \in \delta(i) \backslash j} J_{ik} J_{k \rightarrow i}^{-1} [J_{ki}, h_{k \rightarrow i}],\\
    \mu_i & = J_{i \rightarrow j}^{-1} (h_{i \rightarrow j} - J_{i j} \mu_j),
  \end{aligned}
\end{equation}
where \(J \mu = h\) is the equation that is to be solved. These can be replaced by appropriate square-root forms to obtain instead, an LDL decomposition.\\

\begin{theo}
\label{thm:twlem1}
The linear equation \eqref{eqn:spfamily} can be solved in time \(O(|\mathcal{H}| \tw(\mathcal{H})^3)\), given the tree-decomposition, \emph{via} \Cref{alg:gqp}.
\end{theo}
\emph{Proof.} The correctness of the algorithm follows from \Cref{lem:esep}. The bound holds trivially if, \(\operatorname{rank}{\tilde{G}} \leq \operatorname{rank}{\tilde{Q}_{\iota, \iota}}\), at every step of the algorithm. Otherwise, by realizing that \(\tilde{G}_{l \rightarrow p}\), can't have rank more than \(|\chi(p)|\), the proof follows.

\qed

\end{document}